\RequirePackage{amsthm}
\documentclass[sn-mathphys-num, Numbered]{sn-jnl}

\usepackage{graphicx}%
\usepackage{multirow}%
\usepackage{amsmath,amssymb,amsfonts}%
\usepackage{amsthm}%
\usepackage{mathrsfs}%
\usepackage[title]{appendix}%
\usepackage{xcolor}%
\usepackage{textcomp}%
\usepackage{manyfoot}%
\usepackage{booktabs}%
\usepackage{algorithm}%
\usepackage{algorithmicx}%
\usepackage{algpseudocode}%
\usepackage{listings}%
\usepackage{float}

\theoremstyle{thmstyleone}%
\theoremstyle{thmstyletwo}%

\theoremstyle{thmstylethree}%

\begin{document}

\title[Article Title]{
Constraining Genetic Symbolic Regression via Semantic Backpropagation
}


\author*[1]{\fnm{Maximilian} \sur{Reissmann}}\email{reissmannm@student.unimelb.edu.au}

\author[1]{\fnm{Yuan} \sur{Fang}}\email{yuanfang1@student.unimelb.edu.au}
\author[1]{\fnm{Andrew S. H.} \sur{Ooi}}\email{a.ooi@unimelb.edu.au}
\author[1]{\fnm{Richard D.} \sur{Sandberg}}\email{richard.sandberg@unimelb.edu.au}

\affil*[1]{\orgdiv{Department of Mechanical Engineering}, \orgname{University of Melbourne}, \orgaddress{\street{Grattan Street}, \city{Melbourne}, \postcode{3010}, \state{Victoria}, \country{Australia}}}


\abstract{Evolutionary symbolic regression approaches are powerful tools that can approximate an explicit mapping between input features and observation for various problems. However, ensuring that explored expressions maintain consistency with domain-specific constraints remains a crucial challenge. While neural networks are able to employ additional information like conservation laws to achieve more appropriate and robust approximations, the potential remains unrealized within genetic algorithms. This disparity is rooted in the inherent discrete randomness of recombining and mutating to generate new mapping expressions, making it challenging to maintain and preserve inferred constraints or restrictions in the course of the exploration. To address this limitation, we propose an approach centered on semantic backpropagation incorporated into the Gene Expression Programming (GEP), which integrates domain-specific properties in a vector representation as corrective feedback during the evolutionary process. By creating backward rules akin to algorithmic differentiation and leveraging pre-computed subsolutions, the mechanism allows the enforcement of any constraint within an expression tree by determining the misalignment and propagating desired changes back. To illustrate the effectiveness of constraining GEP through semantic backpropagation, we take the constraint of physical dimension as an example. This framework is applied to discover physical equations from the Feynman lectures. Results have shown not only an increased likelihood of recovering the original equation but also notable robustness in the presence of noisy data.}


\maketitle

\section{Introduction}\label{sec1}
The ability to track, store, and access vast amounts of data, in conjunction with the advent of an increasing landscape of machine learning algorithms, has emerged as a promising avenue to scientific discovery \cite{Brunton2016DiscoveringSystems}, engineering design \cite{LI2022100849}, or biomedical research \cite{Vamathevan2019}. Here, regression approaches not only offer the opportunity to explore underlying patterns but also yield the potential to elaborate data-driven models for enabling promising prediction capabilities. However, a majority of existing methods often fail to achieve generalizability for out-of-training-distribution scenarios or struggle to capture underlying principles that govern natural phenomena when attempting to map the observations \cite{Makke2024}. These issues become even more critical when the number of features increases, implying a vast amount of possibilities to construct such a model, extending the development time and the model uncertainty. 

Strategies such as regularization and constraint-based approaches have been developed to address the generalization gap and narrow the search space. Notably, techniques relying on neural networks are more accessible for advanced regularization due to their internal continuous weight representation and their differentiation mechanism \cite{RAISSI2019686,greydanus2019hamiltonianneuralnetworks}. While these approaches have gained significant traction, there have been only a few attempts within the field of evolutionary computation. Nevertheless, evolutionary approaches, such as Genetic Programming (GP)\cite{koza1992genetic} or Gene Expression Programming (GEP) \cite{9533735}, can handle non-smooth optimization problems and navigate high-dimensional, non-convex error surfaces, making them very versatile and applicable on a wide range of diverse problems. Here, strategies mitigating the effect of an increasing search space can be pointed to methods enforcing constraints \cite{Janikow1996}, such as strict typing \cite{6792634} or grammar-based construction \cite{sobania2021generalizability,ramos2020grammatically}, exist. More recently, success was achieved by employing regularization utilizing domain-specific properties to modify the loss function \cite{bleh2024finding,grundner2024datadriven}.

Still, regularization alone has key limitations, e.g., requiring case-dependent manual tuning of regularization parameters \cite{tikanov} and potentially extending exploration time by suppressing promising candidate solutions in early iterations \cite{bleh2024finding}. On the other hand, methods like strongly typed \cite{stronglytyped} or grammar-guided \cite{grammaroneil} GP necessitate substantial modifications to genetic operators  \cite{mckay} while suffering from low population diversity due to the initialization procedure \cite{inproceedings_grammar_guided_lim}. Moreover, these approaches demonstrate weaker performance compared to other methods \cite{lacava2021contemporarysymbolicregressionmethods}.

Strategies aiming to infer sub-solutions with a more heuristic approach, such as semantic backpropagation \cite{Pawlak2015}, exist but are currently only featuring as a search operator on the regression target.

To address this shortcoming, we propose integrating semantic backpropagation within GEP. The study introduces a methodology for formulating and employing constraints derived from the relationship between inputs and outputs, mapped into a continuous space. Based on that formulation, a metric and the backpropagation mechanism can be utilized to infer the desired changes within the candidate solution. A case study by enforcing dimensional homogeneity of physical units is provided to exemplify the formulation of constraints and the application of backpropagation. By conducting a commonly utilized benchmark for symbolic regression (SRBench \cite{lacava2021contemporarysymbolicregressionmethods}) that, in comparison to other approaches, we demonstrate the framework's capability to enforce and maintain a constraint, limiting the number of possible combinations and, therefore, restricting the search space. 

The study begins with a brief overview of related work. In \autoref{section3:methodology}, we describe the methodology of incorporating physical knowledge, explaining the details of the synergy between the semantic backpropagation and the physical constraints, as well as the details of the algorithmic integration. Thereafter, \autoref{results_and_dists} presents the experimental setup, the employed metrics, and the test scenarios performed, followed by a discussion of the results. The last section concludes this work and outlines future work.

\section{Related Work}
This section starts by defining the relevant terminology and further examines the current state of symbolic regression methodologies, focusing on evolutionary methods. Subsequently, we outline methodologies that incorporate domain-specific knowledge into optimization processes. 
\subsection{Symbolic Regression}
Assuming a regression task is given through a data set consisting of input features $\{\textbf{x}_0,\dots,\textbf{x}_n\} = \textbf{X}$ and a target value $\textbf{y}$, where we seek to find a mapping function $\hat{f}$ so that $\hat{f}(\textbf{X})\approx \textbf{y}$. Employing symbolic regression restricts the model formulation $\hat{f}$ to be an algebraic equation composed of various functions $g_k$ from a pre-defined function space $\Omega$, combining a subset of $\textbf{X}$ and various floating-point numbers or constants. Compared to regression methods relying on neural network architectures, the symbolic ones emphasize their advantage according to interpretability and better generalizability \cite{Makke2024}.

However, these methods come with their challenges, as the problem of identifying an appropriate matching function has been proven to be NP-Hard (nondeterministic polynomial time) \cite{virgolinsnphard}. This implies that assuming $P \neq NP$, whereby $P$ signifies the problems solvable in polynomial time, no algorithm exists to solve all instances in polynomial time of this problem.

Currently, the field of symbolic regression is still primarily dominated by meta-heuristic approaches like genetic algorithms \cite{501943,889734,Cranmer2023,stephens2015gplearn,He2022, koza1992genetic, ferreira2006gene, virgolin2021improving, Kasten2023}. Here, the direction for the optimization finds its role model within biological evolution. A randomly initialized start population, representing a set of different equations, is optimized over many generations. An assessment based on an error metric takes place in each generation, preceding the application of various modifications of the candidate solutions. Building upon Koza's research \cite{koza1992genetic}, subsequent studies have implemented different techniques to enhance the efficiency of the exploration. These include specialized encoding schemes, predefined assumptions, more effective genetic search operators \cite{Pawlak2015,10.1007/978-3-031-43421-1_11}, and strategies for parallel processing and optimized hardware utilization \cite{baeta2021tensorgp}. Despite these improvements, the poor scaling to larger problems, low search efficiency, and unpredictable convergence remain significant issues \cite{Kamienny2022}.

Techniques developed in recent years, like the equation learner (EQL) \cite{pmlr-v80-sahoo18a} and the Deep Symbolic Regression (DSR) Framework \cite{Petersen2021}, aim to improve the search efficiency and suggest a more systematic procedure by using a neural network, hence driving the exploration through gradient information. The previously mentioned approach integrates mathematical operators into a classical neural network forward pass and employs gradient descent. In contrast, DSR utilizes reinforcement learning, a risk-seeking policy, to refine the prediction of a more appropriate equation. Notably, while first-order optimization methods demonstrate proficiency in local search regions, they might end up with a sub-optimal solution. In further studies, this challenge has motivated the development of hybrid methods using the properties of meta-heuristics as well as the guided feedback from a gradient, for instance, the neural guided genetic programming \cite{Mundhenk2021} or the reinforcement GEP\cite{9533735}. 

Yet, the problem of excessive convergence time due to the large number of iterations required to infer a sufficient solution that fits the data remains. Here, the advent of large language models shifts the focus from an intensive search process to the training of parameter-intensive neural networks. Particularly, approaches like \cite{Valipour2021, Kamienny2022, Biggio2020, Biggio2021, Vastl2022} demonstrate that learning a large sample set of equations with their corresponding input- and output spaces allows a close inference, at least for the spatial dimension that has been seen during the training. Furthermore, they also try to approximate numerical constants using a second-order optimization. For example, the End-to-End Symbolic regression proposed by Kamienny et al. \cite{Kamienny2022} demonstrated on the Feynman data set included within the SRBench \cite{lacava2021contemporarysymbolicregressionmethods} that for reaching a comparable result of the $\textbf{R}^2$-score, the inference time is around $100$ times smaller compared to the evolutionary GP-Gomea \cite{virgolin2021improving}. Nevertheless, these approaches become less accurate in the next token prediction when the feature space increases, necessitating enhanced decoding methods \cite{Kamienny2023}, which also increases the risk of an extended exploration time and results in further considerations on mitigating or restricting the search space. Moreover, recent transformer-based approaches for symbolic regression face limitations, mainly when dealing with out-of-distribution datasets \cite{Kamienny2023}.

\subsection{Consideration of a Dimensional Constraint}
In case $\textbf{X}$ is recorded by observing a physical system, each feature $x_i$ represents a physical quantity having a dimension $\phi_i$. For instance, when attempting to determine a velocity, a feature $x_1$ might represent a distance with dimension $\phi_1=[\text{m}]$ (meters) while another feature $x_2$ could be time, having dimension $\phi_2=[\text{s}]$ (seconds). Hence, the target can be derived as the composition of $\phi_1$ and $\phi_2$ signifying velocity ($v=x_1/x_2$) with dimension $[\text{m}/\text{s}]$. This extends the initial given target for $\hat{f}(\textbf{X})\approx \textbf{y}$, to also satisfy dimensional homogeneity $\phi_f = \phi_y$. Imposing such an additional constraint reduces the number of valid combinations between functions and features, implying a restriction on the search space. In the subsequent study, we refer to a candidate solution as invalid if $\phi_f \neq \phi_y$.
A common first step in analyzing a data set involves techniques like normalization or non-dimensionalization of each feature. This process standardizes the range of the quantities within the data set and ensures that the contributions to the final result are comparable. Regardless, in a significant portion of scientific equations, the input features are assigned physical dimensions, representing real-world quantities. The dimensions provide insight into how the quantities are linked to the given solution. However, most of the aforementioned approaches emphasize minimizing a numerical loss between the prediction of the generated function and the target values or trying to achieve a less complex equation regarding the function composition. With this idea in mind, some approaches started incorporating the homogeneity of the physical dimensions as a constraint to restrict the feature space and enable a more consistent representation of possible solutions.

The study of Udrescu and Tegmark \cite{udrescu2020ai} and the subsequent enhanced version \cite{Udrescu2020} describes a framework, referred to as AI Feynman, where the dimensional analysis marks the initial step. A mechanism determines new variables before proceeding to subsequent stages, ensuring a proper scaling of the physical quantities (non-dimensionalization). In a further study, Bakarji et al. \cite{Bakarji2022} proposed three methods, all relying on the Buckingham-Pi theorem \cite{buckingham1914physically} - a principle utilized in engineering science facilitating the identification of appropriate combinations of considered quantities for non-dimensionalization. Before employing Brunton et al.'s proposed approach \cite{Brunton2016DiscoveringSystems}, the author implemented the dimensions as a constraint on the optimization, specifically as a bias prescribing the construction of the first layer of a neural network (BuckiNet). Tenachi et al. \cite{Tenachi2023} extended the DSR approach \cite{Petersen2021} to a physical, symbolic optimization to take the physical dimensions of the data set into account. Here, the constraint is given upfront and acts in the way of a grammar-guided construction rule during the sampling. The authors demonstrate that this improves the rate of finding the baseline solution by approximately a factor of three compared to the original DSR. Particularly in the field of evolutionary methods, Matcheva et al. \cite{matchev2021analytical} demonstrated that initial restrictions through a dimensional analysis provide robust analytical models. Grundner et al. \cite{grundner2024datadriven} incorporate the information about dimensionalized quantities utilizing an implementation referred to as PySR \cite{Cranmer2023} to extract more coherent and robust models within the domain of astrophysics. Within that implementation, the dimensional consistency is enforced by adding a penalty on top of the regression loss.
Further examples, particularly utilizing GEP, are proposed by \cite{bleh2024finding} and \cite{Ma_Zhang_Feng_Xing_Wen_2024} conducting function approximation within the domain of fluid mechanics. While the first proposed method by Bleh et al. \cite{bleh2024finding} discards invalid solutions, reducing the exploration capabilities, Ma et al. \cite{Ma_Zhang_Feng_Xing_Wen_2024} expands the search until the best expression meets the dimensional homogeneity, which further increases the time. 

These observations from the literature strongly support the assumption that dimensional information can guide the exploration of metaheuristics and lead to more precise and robust approximations. Consequently, our approach seeks to establish dimensional homogeneity through a heuristic correction. To the authors' knowledge, no prior study has specifically explored the application of semantic backpropagation to physical constraints, and in particular physical dimensions, nor its effects on exploratory processes. Furthermore, the concept has not been adapted to the methodology of the GEP to date. Therefore, we integrate semantic backpropagation into GEP to enforce the physical unit constraints during the generation of physical equations. The performance of this proposed algorithm will be compared to other approaches mentioned above. 

\section{Methodology}
\label{section3:methodology}
This section starts with a brief overview of the GEP approach and the applied genetic operators. While alternative evolutionary methods exist, GEP was chosen for its efficient handling of complex expressions through genotype-phenotype separation, and due to more robust performance compared to GP \cite{ShiraniFaradonbeh2017,FARADONBEH2016254}. By pointing out changes in the underlying representation of the expression (genotype) caused by modifications through the evolutionary process, we can identify issues that need to be addressed to take constraints into account without restricting the exploration capability. 
Subsequently, we introduce a representation, providing a continuous description of feature properties. This representation is derived from the concrete example involving physical dimensions, enabling the application of a distance metric. The last part of the section introduces the concept of semantic backpropagation and proposes the algorithmic changes incorporated into the evolutionary method.

\subsection{Preliminaries} 
\begin{figure}[t]
  \centering
  \includegraphics[scale=0.7]{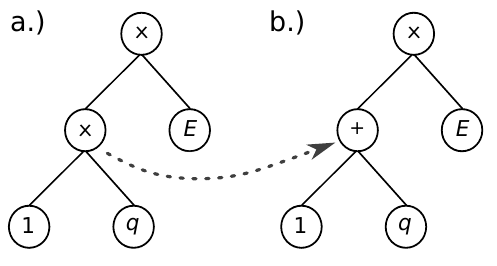}
  \caption{Expression tree built up in pre-order presenting the impact of a mutation operation (indicated by the dashed line), altering the expression from the tree given in a.) to the one represented in b.). Here, the equation sought represents the description of the electric force, given the electric field ($E$) and the electric charge ($q$). Integrating the physical dimension into the evaluation, the resulting tree, particularly the operation $1+q$, leads to a violation.}
  \label{fig1:tree-mutation}
\end{figure}

The methodology of the GEP can be depicted as a particular variant of Genetic Programming \cite{ferreira2006gene, Zhong2017}, whereby the most significant adjustment is characterized through a decoupling of the geno- and phenotype. Here, the former consists of a linear sequence with a fixed size ($2n+1$). The sequence follows the construction rule such that the first $n$-positions, known as the head, may contain both terminal ($x_1, x_2,..$) and non-terminal symbols ($+,-,\times,...$), while the remaining $n+1$, the tail, exclusively includes terminals. We construct the phenotype as an expression tree by decoding the sequence from left to right. For instance, with a head length of two and a genotype of $[\times,\times,1,q,E]$, the phenotype can be determined as illustrated in Fig. \ref{fig1:tree-mutation} a.), whereby the part responsible for constructing the expression is also anointed as Karva expression (K-expression). The K-expression, a linear string, serves as an intermediate representation that bridges the genotype and phenotype, containing only the active symbols that contribute to the expression tree's construction. This linear structure and the decoding mechanism facilitate tracing and the capability of systematic manipulation of the genetic material while preserving structural integrity, which is fundamental to our methodological approach. An overview of the entities employed in the GEP is provided in the appendix (\ref{secA}); for additional details, the reader is referred to \cite{ferreira2006gene}. The optimization is performed by assessing samples of randomly generated expressions (population) and evolving their genotypes by applying different canonical genetic operators, including mutation or crossover. This procedure is repeated over many iterations (generations) until a stop criterion is met. 

While the static definition of the sequence ensures algebraic validity for the usage of the terminals and non-terminals, it does not take into account the other consistency requirements, like the example of the physical dimensional homogeneity. Furthermore, the genetic operators not only facilitate exploration, but they also carry the risk of eliminating viable sub-solutions \cite{Zhong2017}, affecting an underlying system that tries to achieve consistency. 
Referring to the previous example $[\times,\times,1,q,E]$ and assigning the corresponding physical dimension to each symbol, the evaluation of the tree in Fig. \ref{fig1:tree-mutation} a.) yields a correct solution. In detail, $E$ represents the electrical field through $[\text{V}/\text{m}]$ (Volt per meter), $q$ the electric charge given in $[\text{C}]$ (Coulomb) and the value one is dimensionless. Here, the application of the mutation operator leads to a tree with the structure depicted in Fig. \ref{fig1:tree-mutation} b.), which does not align with the expected physical quantity and may be rejected by a constraint despite being possibly close to a ground truth solution. Furthermore, by applying the full range of genetic operators proposed by Fereirra \cite{ferreira2006gene}, the probability of such violations would increase. 

\subsection{Structural Homogeneity}
While the majority of regression techniques often focus on a loss value \cite{ferreira2006gene, Kamienny2022, Kamienny2023, Kasten2023, virgolin2021improving, Valipour2021, Vastl2022}, pointing to the numerical magnitude produced by the mathematical model, formulas, and procedures from areas such as physics \cite{Feynman1963}, chemistry \cite{Reiser2022} or biology \cite{Jha2022} suggest that the structural validity of a model can provide crucial insights. This structure, grounded in scientific principles, not only allows the derivation of the proportion for the contribution of possible input features but also imposes constraints on the possible configurations of the solution. To leverage the information, a mapping is required, satisfying properties such as order and continuity to employ a distance metric, which is able to guide the exploration. Here, a link to a more general approach can be developed by utilizing some key concepts of the expression tree and abstracting it as a graph $G=(V,E)$ (connected graph without cycles) \cite{10.24963/ijcai.2021/595}, where $V$ defines a set of vertices, also referred to as nodes, connected by the edges $E$. To further specify certain properties, techniques such as graph neural networks \cite{10.4108/eetel.3466} assign particular features represented as a continuous feature vector akin to \cite{Reiser2022,Jha2022}. Through defined functions corresponding to the transitional properties of the nodes, the feature vectors can be aggregated. Subsequently, an error metric can be employed to determine the distance between the output and an expected value. Minimizing that discrepancy necessitates a change within the graph. Since we define a constraint, such as the feature vectors are constant according to the source nodes, the change directly affects a node or a segment.  

Throughout various domains, such a feature vector can be developed, for instance, representing the orbital states of the valence configuration when dealing with chemical reactions \cite{Reiser2022} or language model extracted protein sequences within a biomedical area \cite{Jha2022}. For the example investigated within this study, we demonstrate the feature mapping as well as the definition of the function utilized for the aggregation on the physical dimension, particularly the international system of the physical basis units (SI-unit System). 

This approach provides a unique system for standardizing measurements and making their values comparable across different domains \cite{Wiersma2021}. It consists of seven fundamental units, each representing a specific physical dimension: mass (kilogram - $\text{[kg]}$), length  (meter - $\text{[m]}$), time (second - $\text{[s]}$), temperature (Kelvin - $\text{[K]}$), current (Ampere - $\text{[A]}$), quantity (mole - $\text{[mol]}$), and luminous intensity (candela - $\text{[cd]}$). Their composition via mathematical operators leads to more complex representations of a physical quantity. According to the application for existing optimization frameworks, packages available can be found within the PySR implementation \cite{Cranmer2023} or provided through SciPy \cite{2020SciPy-NMeth}. 
\begin{figure}[t]
  \centering
  \includegraphics[scale=0.3]{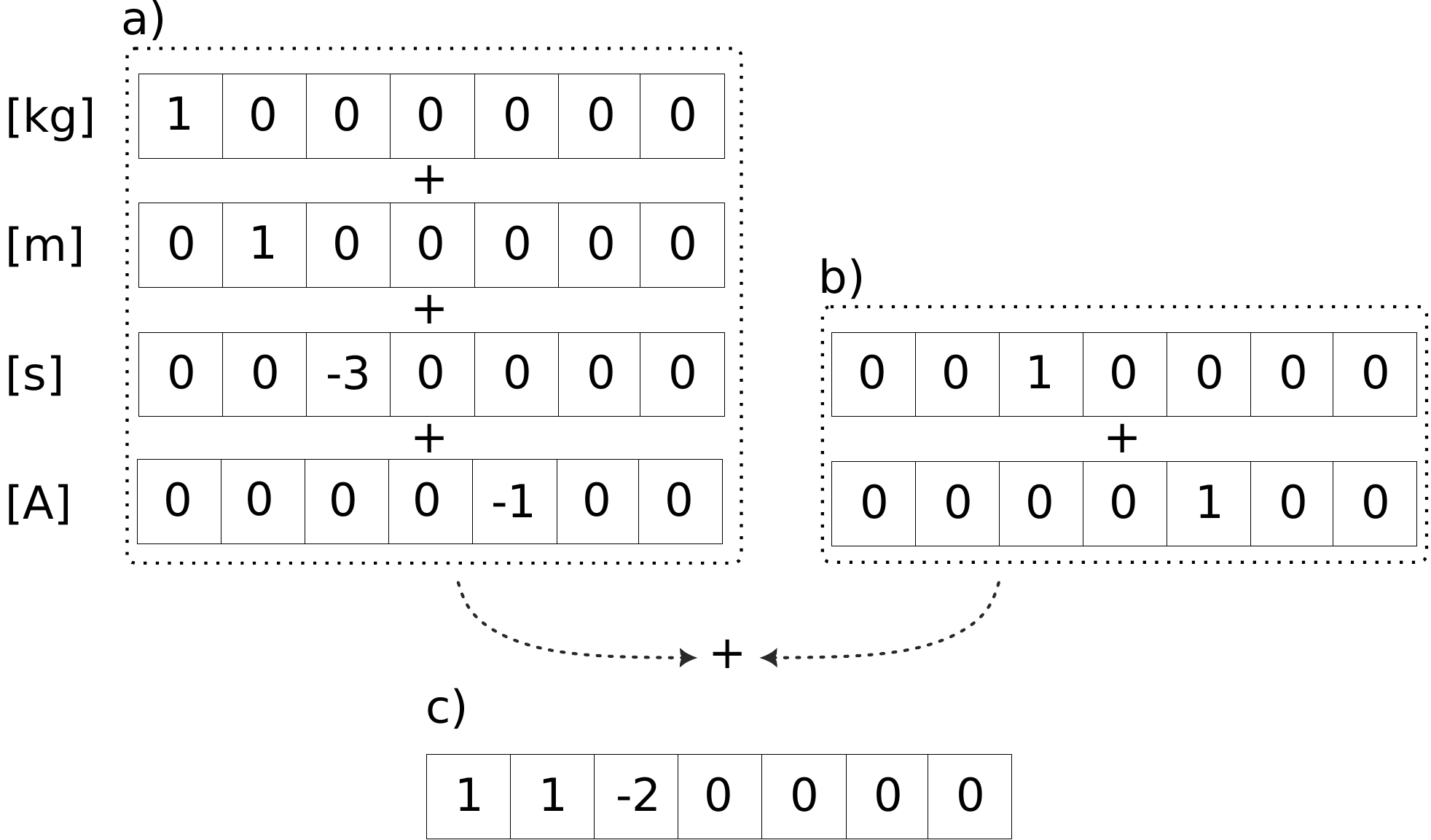}
  \caption{Representation and composition of different physical dimensions based on the basis units. Here, a.) depicts the dimension of the electrical field $E$ and b.) the electrical charge $q$. The employed operation on the two quantities results in a force, which can be described by the composition in c.).}
  \label{fig:simple_si_units}
\end{figure}
Similar to Bleh et al. \cite{bleh2024finding}, we utilize a system based on a vector representation. This is not only due to its simplicity but also because it provides an essential part of establishing a distance metric concerning the semantic backpropagation methodology. As done in the software OpenFOAM \cite{OpenFOAM}, the information is projected into a one-hot encoded vector with a length of seven, whereby each position corresponds to one of these dimensions. 

Based on that representation, more complex relations can be created by increasing or decreasing the number at each position. Referring to the example given in Fig \ref{fig1:tree-mutation} a.), where Coulomb's law is expressed as the product of the electric charge ($q$) and the electric field ($E$), the dimension can be calculated as demonstrated in Fig. \ref{fig:simple_si_units}. Subfigure a.) signifies the dimension for $E$. The operations for the dimensions related to the multiplication of the features lead to the addition yielding the corresponding unit ([$\text{kg}\cdot\text{m}\cdot\text{s}^{-3}\cdot\text{A}^{-1}$] $\implies$ [$\text{V/m}$]). A comprehensive explanation of various SI-unit compositions is provided in Appendix \ref{secC}. Subfigure b.) depicts the representation for $q$, obtained from ([$\text{s}\cdot\text{A}$] $\implies$ [$\text{C}$]). Finally, Subfigure c.) illustrates how the mathematical operation of the two quantities induces a force, which is represented by its physical dimension ([$\text{kg}\cdot\text{m}\cdot\text{s}^{-2}$] $\implies$ [$\text{N}$]). 

Given that the mathematical operator on the particular symbol does not satisfy the evaluation of the physical dimensions, we define a set of rules for performing the assessment. Here, Table \ref{tab:symsops} lists the transformation for the standard functions employed within the utilized framework. The letters \textbf{a} and \textbf{b} signify two arbitrary physical quantities with their corresponding $[\cdot]$ dimensional representation. Applying addition or subtraction requires that both values correspond to the same physical dimension. We can employ arithmetic functions as depicted in Table \ref{tab:symsops} for the additional operators. Particularly, the forward operations are related to the evaluation path from the input to the output, and the backward one describes the utilization from the output to the inputs.
\begin{table}
  \caption{Physical dimension transformation rules for evaluating expression in both forward and backward directions. $\phi_c$ represents the resulting physical dimension, while $\phi_a$ and $\phi_b$ mark the inputs.}
  \label{tab:symsops}
  \begin{tabular}{llll}
    \toprule
    Operator & Forward & Backward  & Comment\\
    &Rule& Rule&\\
    \midrule
    $c=a + b$ &  $\phi_c:= \phi_a = \phi_b$& $\phi_a:= \phi_c $ & Entries of a and b needs to be\\
    &&$\phi_b:= \phi_c $ & the same on every position \\
    $c=a - b$ &  $\phi_c:= \phi_a  = \phi_b$& $\phi_a:= \phi_c $  & ref. '+' operation\\
    &&$\phi_b:= \phi_c $ & the same on every position \\
    $c=a * b$ & $\phi_c:= \phi_a + \phi_b$& ref. appendix & Pointwise addition\\
    $c=a / b$ & $\phi_c:= \phi_a - \phi_b$& ref. appendix & Pointwise subtraction\\
    $c=\text{a}^n$ & $\phi_c:= n * \phi_a$& $\phi_a:=\phi_c/n$ & Dimension is multiplied by $n$\\
    $c=\sqrt{a}$ & $\phi_c:= \phi_a / 2$ & $\phi_a:=2*\phi_c$ & Dimension is divided by two\\
    $c=sin(a); c=cos(a)$; .. & $\phi_c:= 0*\phi_a $ & $\phi_a:=0*\phi_c$ & Result becomes dimensionless\\&&&Input needs to be dimensionless\\
    $c=log(a); c=exp(a)$; .. & $\phi_c:= 0*\phi_a $ & $\phi_a:=0*\phi_c$ & Result is dimensionless\\&&&Input needs to be dimensionless\\
  \bottomrule
\end{tabular}
\end{table}

\subsection{Semantic Backpropagation}
As indicated in Fig. \ref{fig1:tree-mutation}, the application of genetic operators needs to be encapsulated within a mechanism that verifies dimensional homogeneity and corrects any detected discrepancies. Besides those mentioned above and grammar-guided approaches like \cite{942529, Tenachi2023, Sun2019}, methods that take intermediate states and program semantics into account exist. Here, semantic backpropagation, introduced by Pawlak et al. \cite{Pawlak2015} and applied to the field of genetic programming \cite{virgolin2021improving}, provides a technique on how to manipulate the expression according to a particular purpose, referred to the desired semantics. Expressions are typically evaluated using a hierarchical computational graph depicted as a tree, which can be decomposed into sub-trees $\hat{f}(x)=\hat{f}_2(\hat{f}_1(x),x)$. 
Traversing from the leaves upward, each node contributes an intermediate result, enabling the tracing of the performed executions and relating them to the output. Sub-trees yielding identical results demonstrate shared semantics, which can simplify expressions without losing their semantic value. Moreover, semantic backpropagation assumes that the desired semantics for a function $\hat{f}$ can be achieved by retaining the structure of $\hat{f}$ and adapting iteratively the inner functions $\hat{f}_2$ and $\hat{f}_1$ through heuristic inversion \cite{forwardProp}. A lookup into a library can be performed at each node to determine whether the sub-tree can be replaced by a composition meeting the desired semantics, similar to the rule application for algorithmic differentiation \cite{Wang2023}. The library can be established either in advance, based on the features and a maximum tree height, or dynamically within the course of the exploration. 

\begin{figure}[t]
  \centering
  \includegraphics[scale=.75]{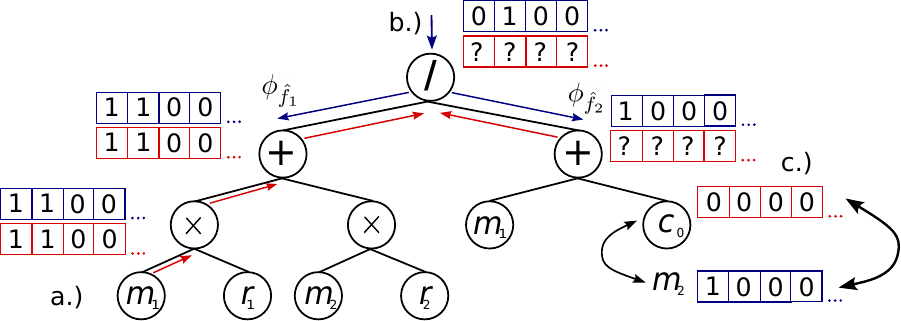}
  \caption{Expression tree built up in pre-order presenting an example approximation for exploring the equation that calculates the center of gravity. The blue boxes illustrate the vector presentations of the desired physical dimension, while the red boxes present the calculated units. Starting from a.) and following the red arrows towards the output, the dimension can be evaluated. A misalignment b.) initiates a process, propagating residual information to the underlying sub-trees (blue arrows), where changes are enforced, here at c.). 
  }
  \label{fig2:tree-mutation}
\end{figure}

Fig. \ref{fig2:tree-mutation} depicts an approximation $\hat{f}$ close to the equation that calculates the center of gravity depending on two different masses ($m_1,m_2$) and their corresponding position ($r_1,r_2$). The desired semantics ($\phi_{\hat{f}}$), here provided through the vector representation of the physical dimension, are given with the problem that requires to be solved. While the expected output ought to be a `length', the evaluation by executing the forward path (red arrows) leads to an invalid dimension (red box at b.)). Instead of discarding the whole expression like in \cite{bleh2024finding}, that information, according to the misalignment, is utilized to modify the expression. After quantifying the distance, the backpropagation process is initiated (b.)). The first node encountered signifies a divide operator that splits $\hat{f}$ into $\hat{f}_1$ and $\hat{f}_2$. In the subsequent layer, the intermediate states $\phi_{\hat{f}_1}$ and $\phi_{\hat{f}_2}$ are analyzed along with the algebraic rules defined in Table \ref{tab:symsops}, whereby backward signifies the direction from the output to the input. As $\phi_{\hat{f}_2}$ is undefined, caused by a violation of the desired constraint, here an underlying mismatch of the dimension, the residual for further propagation is calculated by: $-[?,?,?,?,\dots]=[0,1,0,0\dots] - [1,1,0,0,\dots]$, resulting in $[1,0,0,0,\dots]$ as desired semantic for $\phi_{\hat{f}_2}$. This procedure is repeated in the subsequent subtree, where the misalignment can be pointed to the application of the addition function to $m_0$ ($[1,0,0,0,\dots]$) and an arbitrary constant $c_0$ ($[0,0,0,0,\dots]$), resulting in a substitution of an equivalent symbol.
The example outlines the particular algorithmic traits necessary to incorporate the concept into GEP. These mechanisms encompass establishing a library containing different semantics, selecting an appropriate distance metric, identifying the genetic entity to which the corrections are applied, and establishing rules to guide information propagation to subsequent branches:
\begin{enumerate} 
\item Library of Semantics: The library can be described mathematically as a multiset partition ($\mathbf{S}$) of the union of the terminal and the non-terminal symbol set, where each element ($\mathbf{S}_i \in \mathbf{S}$) satisfies defined conditions. These include that $\mathbf{S}_i$ can be transformed into a constant or an algebraic equation, yielding a valid physical dimension $\phi$, and does meet the restriction according to the element count ($\forall i\in \{1,\dots,|S|\}, |\mathbf{S}_i|\leq$ head length), analogous to the maximum tree depth in GP \cite{Pawlak2015,virgolin2021improving}. Utilizing $\phi$ and $|\mathbf{S}_i|$, this can transform $\mathbf{S}$ into a hash table structure whereby the key, defined by these quantities, is assigned to a list containing multisets ($key_{\phi,|\mathbf{S}_i|} \rightarrow {\hat{f}_0,\dots,\hat{f}_i}$). During the runtime, the substitution procedure fetches an entry based on the key and randomly selects an entry from the obtained list.
\item Distance Metric: This measure indicates how close the semantics of $\hat{f}$ and $f$ are, and it further provides a value to fetch the closest candidate from the library. Within this study, the measure $d_\phi$ applied for determining the distance between two vector representations $\phi_a$ and $\phi_b$ is given by the Euclidean distance (elementwise):
\begin{align}
    d_\phi(\phi_a, \phi_b) = \frac{1}{7}\sum\left((\phi_a - \phi_b)^2 \right) \,.
\end{align} 
\item Adaptable entities: A core aspect of GEP's internal mechanisms is that all genetic operators are employed on the genotype representation only. To ensure dimensional homogeneity through phenotype evaluation, changes inferred via semantic backpropagation are applied to the K-expression. The overall gene length constrains the maximum length of any particular replacement. While GP relies on techniques like `Random Disered Operator' (RDO) or `Approximately Geometric Semantic Crossover' (AGX) \cite{rdo_krawiek}, which randomly selects a node, the introduced operator follows the expression tree in a systematic traversing order. Specifically, it determines the potential replacement of a sub-tree according to the desired semantics, leveraging the existing database.  
\item Backpropagation rules: Traversing through the expression, at each node, a decision is required on how to split the information to propagate further or enforce a change based on a matching library entry. Here, some ideas derived from algorithmic differentiation are utilized \cite{Wang2023}. Next to the standard forward execution, we define a backward pass, enabling the capability to delegate residuals to the intermediate states. This residual contains a fraction of the backpropagated vector and is calculated by a defined backward rule (ref. Table \ref{tab:symsops}).
\end{enumerate}
Due to the inherent nature of heuristic approaches, the correction process may not achieve the expected dimensional homogeneity. This scenario could arise if, for instance, the problem is not fully described in terms of the features, the library does not contain a suitable candidate, or the parameter selected for the GEP is insufficient to create an equation appropriate for the problem. Apart from an adjustment of the hyperparameters, a regularization in the form of a loss penalty similar to \cite{grundner2024datadriven,Bakarji2022} or the formulation of a further objective can be introduced to guide the exploration:
\begin{align}
\label{eqn:regularization}
    \hat{\mathcal{L}}(\mathbf{Y}, \mathbf{X}, \hat{f}) = \left(  \mathcal{L}(\mathbf{Y}, \hat{f}(\mathbf{X})) + \lambda \cdot d_{\phi}(\phi_{\mathbf{Y}}, \phi_{\hat{f}(\mathbf{X})})  \right)\, .
\end{align}
Here, $\mathcal{L}$ signifies an arbitrary loss function, $\hat{f}$ an equation approximated by the GEP, and the variable $\lambda$ controls the magnitude of the penalty term, defining the magnitude of how the dimensional homogeneity influences the search.

The subsequent algorithm (ref. \ref{algo:semantic}) summarizes the detailed steps undertaken to preserve dimensional homogeneity. Here, we utilize the pre-computed library ($\mathbf{L}$) to perform the correction for a population ($\mathbf{P}$). The objective is to ensure that the corrected candidates show a high probability of matching the physical dimension of the target $f$, within a small tolerance $\epsilon$. The `propagate\_change' function recursively attempts to correct sub-trees according to the inserted dimension. The function initializes by setting the homogeneity variable to false, indicating that the dimensional consistency condition is not yet satisfied. This is followed by evaluating whether the dimension of the inserted expression $\phi_g$ already satisfies the target dimension $\phi_t$ (Line 3 - Algorithm \ref{algo:semantic}). If this condition is met, the function terminates and returns true, bypassing subsequent `else' branches. Otherwise, the algorithm proceeds by retrieving the left child node of g ($\text{g}_l$) and computing its corresponding target dimension ($\phi_{tl}$) based on the arithmetic operation rule associated with g (Line 7 - Algorithm \ref{algo:semantic}). When the dimension of ($\text{g}_l$) deviates from $\phi_{tl}$, the algorithm attempts to fetch an appropriate replacement from the library $\mathbf{L}$ (Line 12 - Algorithm \ref{algo:semantic}) before propagating the dimensional correction to the underlying child nodes of ($\text{g}_l$). Upon reaching a termination criterion for ($\text{g}_l$), the function applies an identical sequence of operations to the right child node of g. To mitigate the disruptive effects of genetic operators on dimensional homogeneity, the correction process is implemented before evaluation takes place, using several attempts for the correction.
\begin{algorithm}[H]
\caption{Perform Semantic Backpropagation for Dimensional Correction}
\label{algo:semantic}
\begin{algorithmic}[1]
\Require $\mathbf{P}, \, \mathbf{L}$ 
\Ensure $d_\phi(\phi_{\hat{f}_i},\phi_f) < \epsilon, \quad \forall \hat{f}_i \in \mathbf{P}$
\Function{propagate\_change}{$\text{g}, \phi_{t}$}
    \State homogeneity $\Leftarrow$ false
    \If{$d_\phi(\phi_{t}, \phi_{\text{g}}) < \epsilon$}
        \State homogeneity $\Leftarrow$ true \Comment{Set true if dimensional homogeneity is satisfied} \
    \Else
    \State $\text{g}_l \Leftarrow get\_node\_left(\text{g})$
    \State $\phi_{tl} \Leftarrow \text{rule}(\text{g})$ \Comment{determine node target by rule - \ref{tab:symsops}}
    
    \If{$d_\phi(\phi_{tl}, \phi_{\text{g}_{l}}) < \epsilon$} 
    \State homogeneity $\Leftarrow$ true \Comment{homogeneity for left node satisfied} \
    \Else
    \If{$\phi_{tl} \in \mathbf{L}$} \Comment{Appropriate replacement in library?}
    \State $\text{g}_{l} \Leftarrow get(\mathbf{L},\phi_{tl} )$ \Comment{retrieve replacement for node}
    \Else
    \State homogeneity $\Leftarrow$ propagate\_change({$\text{g}_{l}, \phi_{tl}$})
    \EndIf
    \EndIf

    \State $\text{g}_r \Leftarrow get\_node\_right(\text{g})$
    \State $\phi_{tr} \Leftarrow \text{rule}(\text{g})$ \Comment{determine node target by rule - \ref{tab:symsops}} 
    
\If{$d_\phi(\phi_{tr}, \phi_{\text{g}_{r}}) > \epsilon$} 
    \If{$\phi_{tr} \in \mathbf{L}$}
        \State $\text{g}_{r} \Leftarrow get(\mathbf{L},\phi_{tr} )$ \Comment{retrieve replacement for node}
    \Else
        \State homogeneity $\Leftarrow$ propagate\_change({$\text{g}_{r}, \phi_{tr}$}) $\land$ homogeneity
    \EndIf
\EndIf
    \EndIf
    \State \Return homogeneity
\EndFunction

\For{each $\hat{f}_i \in  \mathbf{P}$}
    \State $\text{g}_i \Leftarrow \text{convert}(\hat{f}_i)$ \Comment{create the computational graph}
    \For{each $\_ \in \text{cycles}$} \Comment{pre-defined number of cycles}
        \If{\Call{propagate\_change}{$\text{g}_i, \phi_f$}}
            \State break
        \EndIf
    \EndFor
\EndFor
\end{algorithmic}
\end{algorithm}
\newpage
\section{Results and Discussion}
\label{results_and_dists}
In the subsequent section, we quantify and discuss the influence of dimensional homogeneity on the performance of genetic exploration. We provide a brief overview of the employed benchmark, the adapted parametrization of the GEP implementation, and the utilized error metrics. Subsequently, we first compare the influence of the proposed extension on the GEP within a small test environment to evaluate its effects by focusing on error distribution and functional complexity. Here, the term small refers to a limitation of epochs and the population size. In a subsequent consideration, we then employ the configuration yielding the highest performance and compare it against different approaches from the literature using a more competitive amount of candidate solutions and a larger population size. The hardware utilized was a 32-core server (CPU-AMD-3.6GHz). 

\subsection{Benchmark}
To evaluate the performance of the proposed extension, we utilize a benchmark for symbolic regression established by La Cava et al. \cite{lacava2021contemporarysymbolicregressionmethods}, referred to as SRBench. The benchmark provides the results of various state-of-the-art methods such as GP-Gomea \cite{virgolin2021improving}, AI-Feynman, or the DSR approach conducted on different datasets. Table 2 lists the approaches from the set. The repository comprises 252 datasets with equations from the Feynman Symbolic Regression Database (Feynman Lectures on Physics \cite{Feynman1963}), the Strogaz repository, and several real-world data sets outlined as Black-Box-Problems. The datasets are categorized into different difficulty levels, referred to as `easy,' `medium', or `hard'. Additionally, the SRBench suggests a set of rules and metrics according to the comparability. The incorporated rules are an upper limit of evaluations ($10^6$), the maximum time, the number of independent trials, and the consideration of the train and test split for the particular cases. 
Following the benchmark guidelines, we conducted approximately 100 instances from the Feynman dataset across ten independent trials. For each trial, we randomly shuffled and split a total of 10,000 data points into training and test sets, using a 0.75/0.25 ratio, respectively. The maximum number of candidate evaluations during the exploration is set to the mentioned upper bound. To facilitate dimensional analysis, the benchmark cases were extended with a dictionary providing the dimension information for each symbol within each equation.

We also consider different levels of noise sampled from a normal distribution utilizing the root mean square value and add it to the target \cite{lacava2021contemporarysymbolicregressionmethods}:
\begin{align}
\label{eqn:noise}
\Tilde{y_i} = y_i + \eta \, | \, \eta \sim \mathcal{N} \left ( 0, \gamma\sqrt{\frac{1}{N}\sum y^2} \right ) \, .
\end{align}
Here, $\gamma$ takes the values $\{0.001, 0.01, 0.1\}$, evaluating how well the method performs under different levels of uncertainty. 

The assessment takes place by using the suggested metrics, which includes the $\textbf{R}^2$-score (elaborated on the validation set):
\begin{align}
\label{eqn:r2-score}
R^2(\textbf{Y}, \textbf{X} ,\hat{f}) = 1 - \frac{\sum_{i=1}^{n} (y_i - \hat{f}(x_i))^2}{\sum_{i=1}^{n} (y_i - \bar{y})^2} \, .
\end{align}
This metric describes how well the model predictions explain a portion of the variance in the dependent variables of the data set, whereby a value of one indicates a perfect fit, and zero or less suggests inadequate explanatory capabilities. Additionally, performance measurements include the run time, the average complexity of the approximation, and the accuracy according to the recovery of the ground truth symbolic expression. The first mentioned is determined by counting the contained symbols within the simplified equation, whereby we employ the `Symbolics.jl' library \cite{gowda2021high} instead of `SymPy' \cite{10.7717/peerj-cs.103}. In terms of accuracy, given an approximation $\hat{f}(\mathbf{x},\hat{\theta})$ expressed as functional compositions using the features $x_n$ and the explored coefficients $\hat{\theta}$, symbolic accuracy ($||f-\hat{f}||_\Psi$) can be calculated as follows \cite{lacava2021contemporarysymbolicregressionmethods}:
\begin{align}
\label{function:solacc}
    ||f-\hat{f}||_\Psi = \begin{cases} 
    \text{True}, & f-\hat{f} = 0 \lor f-\hat{f} = a, \text{ where } a \in \mathbb{R} \text{ is constant}, \\
    \text{True}, & f/\hat{f} = 1 \lor f/\hat{f} = a, \text{ where } a \in \mathbb{R} \text{ is constant}, \\
    \text{False}, & \text{otherwise}. \\
\end{cases}
\end{align}
Averaging the occurrence of $||f-\hat{f}||_\Psi$ over the number of tests then expresses the rate of how often the algorithm could find the true equation (symbolic solution rate). 

\subsection{Hyperparameter Setup}
The trials are performed utilizing an in-house version of the GEP (implemented in Julia \cite{bezanson2012julia}), with a parameter configuration derived from previous suggestions \cite{ferreira2006gene} and outlined in Table \ref{tab:high}. It is worth noting that both the solution quality and the convergence time are highly sensitive to these hyperparameter settings, which ideally are tailored to each problem's specific complexity. In our research, the primary objective is to examine the effects of the consideration of physical dimensions rather than to surpass the performance of existing algorithms.
\begin{table}
  \caption{Parameterization of the Evolutionary Method. The first column delineates the parameter name, the second column yields the corresponding value, and the third column provides additional context.}
  \label{tab:high}
  \begin{tabular}{lll}
    \toprule
    Parameter & Value & Comment\\
    \midrule
    Population size & 500 / 1000 & Maximum number of\\
    (GEP\_Test / SRBench) & & candidate solutions\\
    Generations & 1000 / 1500 & Maximum number of\\
    (GEP\_Test / SRBench) & & iterations\\
    Head length & 8 & Maximum positions \\
    & & non-terminals in a gene\\
    Number of genes & 3 & Number of sub-expressions\\
    Tournament size & 3 & Number of contenders in the\\
    & & tournament\\
    Mating proportion & 0.5 & Fraction of the population for\\
    & & for recreation\\
    Mutation probability & 0.2 & Chance of random gene\\
    & & modification\\
    Inversion probability & 0.1 & Probability of sequence\\
    & & reversal (classic GEP operator)\\
    One-point crossover & 0.5 & Chance of single-point\\
    probability & & genetic exchange\\
    Two-point crossover & 0.4 & Chance of two-point\\
    probability & & genetic exchange\\
    Non-terminal set & $\{+,-,*,/,\log,\exp,\sin(),$ & Available mathematical\\
    & $\cos(),()^2,\sqrt{()}\}$ & operators\\
    Terminal set & $x_1,...,x_n$ & Available input features\\
    Coefficient-Optimizer & Conjugated gradient & Method for constant\\
    & & optimization\\ 
    \bottomrule
  \end{tabular}
\end{table}
The fitness measurement of a candidate solution to drive the exploration is given by mean squared error (MSE) with an adaptable regularization term:
\begin{align}
\hat{\mathcal{L}}(\textbf{Y},\textbf{X},\hat{f}) = \frac{1}{n} \sum_{i=1}^{n} (y_i - \hat{f}(x_i))^2 + \lambda \cdot d_\phi\,.
\label{eqn:rms}
\end{align}
Since the comparison also includes the widely applied approach of discarding or penalizing the fitness value, the error function is adapted using the following regularization:
\begin{align}
\label{function:dimension}
    d_\phi(\textbf{Y},\textbf{X},\hat{f}) = \begin{cases} 
\|\phi_{\textbf{Y}}-\phi_{\hat{f}(\textbf{X})} \|_2, & \text{homogeneity violated}, \\
\infty, & \text{otherwise}.
\end{cases}
\end{align}
The first line of the equation (\ref{function:dimension}) evaluates an L2-norm between the resulting dimension evaluating $\hat{f}$ and the dimension of the target values $\textbf{Y}$. The choice of the $\lambda$ values, taken from the set $\{0,0.1,1,10\}$, is intended to study the effect of the search direction. A value of $0$ for $\lambda$ represents the baseline method. The intermediate values are likely to influence the metric in smaller error ranges, while the higher values could dominate even in more extensive error regimes (when $\hat{\mathcal{L}}\approx10^3$). The second line is applied when the evaluation is invalid due to a mismatch of the utilized algebraic operations and the dimensions similar to Fig. \ref{fig2:tree-mutation}.
\subsection{Evaluation}
\begin{figure}[t]
  \centering
  \hspace*{0.65cm}\includegraphics[scale=0.33]
  {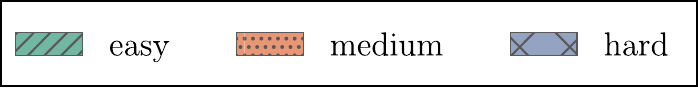}\\
  \vspace{0.2cm}
  \centering
  \hspace*{-0.2cm}\includegraphics[scale=0.27]{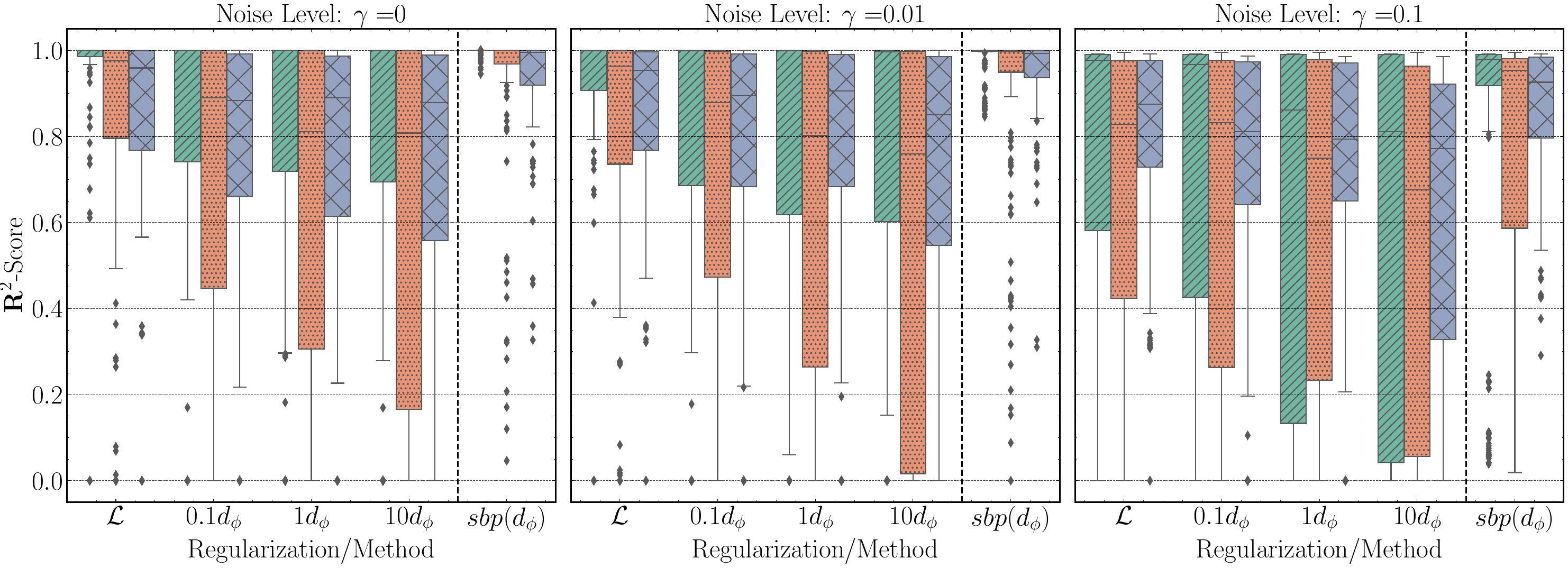}
  \caption{Boxplot, $\textbf{R}^2$-score comparison of the best-explored candidate solutions (y-axis) between the proposed (SBP) and standard method GEP with varying fitness regularization ($\lambda \in \{0.1,1,10\}$) (x-axis) across increasing levels of noise $\gamma \in \{0, 0.01, 0.1\}$. The individual colors for each category signify the separation into a degree of difficulty (green-/), `medium' (orange-..), and `hard' (blue-x). Here, a score of one outlines a perfect match.}
  \label{plot1:r2-test}
\end{figure}
Before comparing the proposed approach with existing methods in the literature, we first examine the capabilities of our introduced technique with respect to the standard GEP and the different quantities of the regularization. 
\subsubsection{Comparative Analysis of GEP}
Analyzing the $\textbf{R}^2$-scores calculated from the test dataset across the conducted test trials and further separating them into their distinct categories (Fig. \ref{eqn:r2-score}), we observe that candidate solutions employing fitness regularization underperform compared to both the standard method ($\mathcal{L}$) and the proposed method ($sbp(\phi)$). Contrary to the expectation based on the difficulty level, this strategy ($\lambda \in \{0.1,1,10\}$) reveals its highest dispersion for the category labeled as `medium' (orange bar), particularly for $\lambda=10$. Within the outlined test scenario for $\gamma=0$, the median for the mentioned category shifts from $0.87$ to $0.74$, up to $22\%$ lower than the standard GEP. The performance further deteriorates when $\gamma>0$ is inferred onto the target value, also considerably affecting the `easy' samples while only slightly degrading the `hard' ones, on average from $1.0$ to $0.86$ and $0.88$ to $0.78$, respectively. 

This phenomenon relates to the magnitude of the regression error in combination with the violation that occurred. Since the `easy' samples often consist of only two to three features, the chance to achieve dimensional homogeneity is relatively high, even without correcting them, making them less likely to be dominated by loss regularization. The `hard' samples, on the other hand, face a high probability of violating that property. Still, due to their complex nature and the difficulty of exploring them, the regression loss exceeds the magnitude of the violation penalty, leading to a negligible influence. In contrast, most of the `medium' samples also yield a considerable degree of freedom. Yet, the likelihood is higher that the approximation shows a loss where the magnitude is close to the regularization, indicating a high efficiency of the search algorithm. By following this trend, we can anticipate that further incrementing ($\lambda \to \infty$), which emulates the discarding of a candidate solution akin to \cite{bleh2024finding}, mitigates that success according to that score. This results in an extension of the exploration time to achieve a similar result to the standard approach. This observation aligns with the literature, where regularization increases training efforts.  

Comparing the proposed technique, which aims to ensure dimensional homogeneity, to the standard method, the median of the score not only is increased for the varying difficulty (easy, medium, hard) cases but also outperforms the standard procedure on the three depicted noise levels. This is reflected in a decrease of the dispersion, a shift of the median towards a higher $\mathbf{R}^2$-score, and a mitigation of the negative influence of the noise. Notably, for the most challenging samples (purple bar), our novel approach improves the performance by $4.24\%$, $4.87\%$, and $13.45\%$ for $\gamma \in \{0,0.01,0.1\}$, respectively. 

\begin{figure}[t]
  \centering
  \hspace*{0.65cm}\includegraphics[scale=0.33]
  {images/plot_1_legend.pdf}\\
  \vspace{0.2cm}
  \centering
  \hspace*{-0.2cm}\includegraphics[scale=0.27]{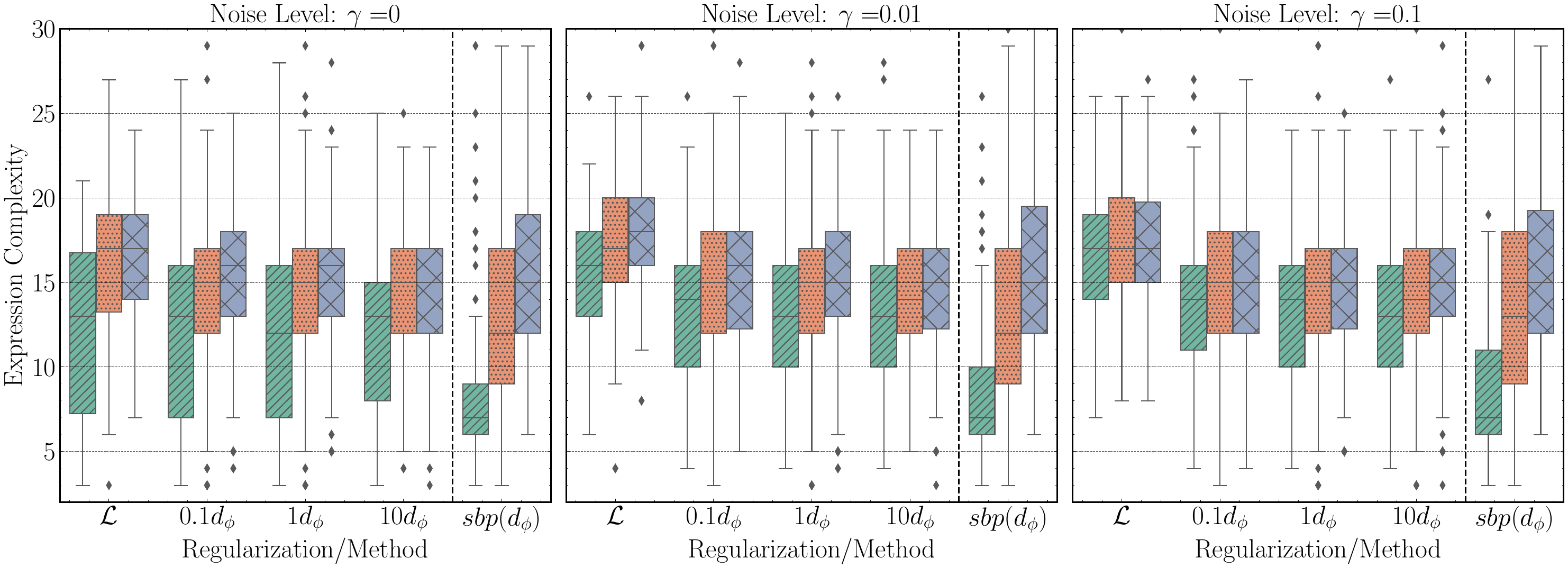}
  \caption{Comparison of the complexity of the expression for the final generation using the specified function from the \textbf{Symbolics.jl} library \cite{gowda2021high} (y-axis).  
  Here, the expression complexity is further outlined for different levels of added noise ($\gamma \in \{0, 0.01, 0.1\}$), the varying regularization ($\lambda \in \{0.1,1,10\}$), and the introduced method. The individual colors for each category signify the separation into a degree of difficulty: `easy' (green-/), `medium' (orange-..), and `hard' (blue-x). A smaller number indicates a less complex expression.}
  \label{plot2:complexity}
\end{figure}

Model complexity, particularly size, serves as an additional criterion for assessing solution quality in symbolic regression. Beyond accuracy comparisons with ground-truth solutions,  complementary studies like \cite{Ma_Zhang_Feng_Xing_Wen_2024, bleh2024finding, grundner2024datadriven} indicate that the size of the model relates to the robustness. Larger models are more prone to overfitting the given dataset, are less reliable when extrapolating beyond the range of the training data, and are more difficult to interpret. Here, Fig. \ref{plot2:complexity} offers insight into the complexity of the explored models achieving the scores mentioned before. The comparison accounts for the variation of the chosen penalty and the increasing noise level. As expected, independent of the adaptation, the pattern of the median complexity shows coherence with the difficulty level, marking the lowest scores for the `easy' ones and higher values for the `hard' examples. Introducing noise ($\gamma=0.1$) induces subtle shifts in the respective medians towards more complex functional representations when only relying on the data as a source of information. This effect is particularly noticeable when employing the standard method for 'easy' equations, yielding complexity ranges of $14.0-18.0$. In contrast, methods utilizing loss penalties or dimensional correction encounter only slight changes of $8\%$ or remain stable at seven. 

\begin{figure}[t]
  \centering
  \hspace*{0.65cm}\includegraphics[scale=0.6]{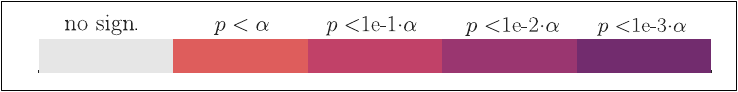}\\
  \vspace{0.4cm}
  \centering
  \hspace*{-0.2cm}\includegraphics[scale=0.6]{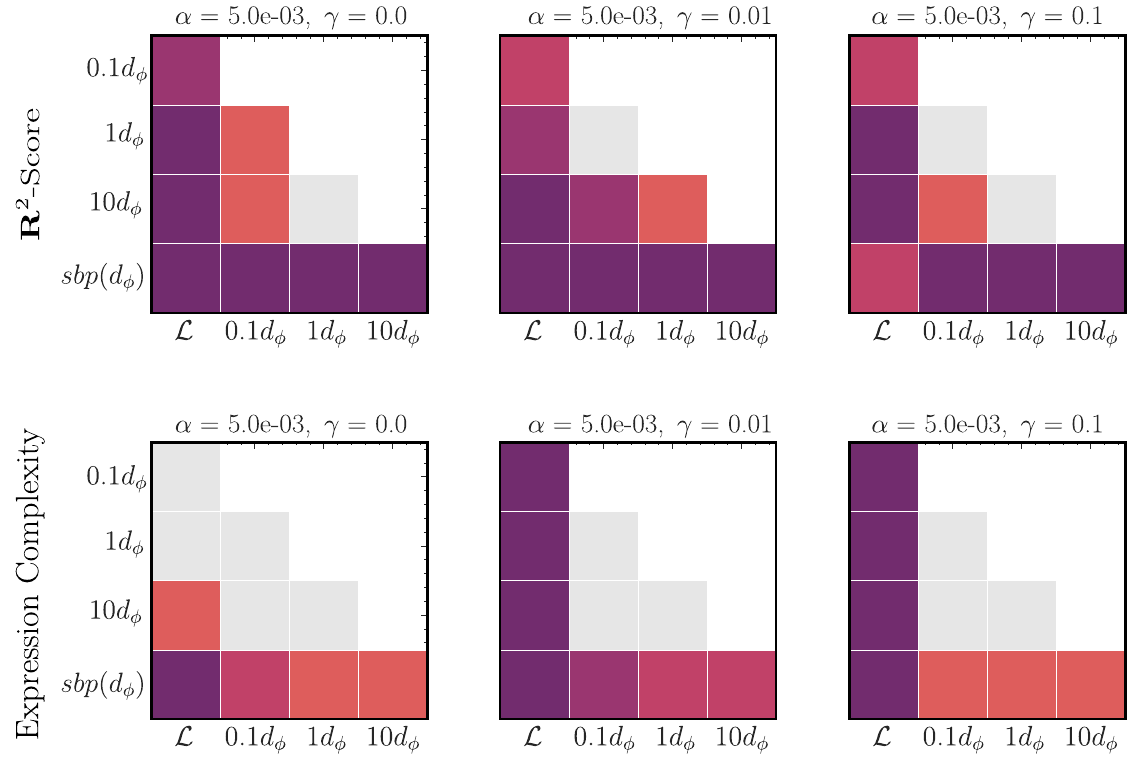}
  \caption{Heatmap outlining the pair-wise comparison according to the statistical significance over the different levels of target noise ($\gamma \in \{0,0.01,0.1\}$) for the $\mathbf{R}^2$-Score (plots first row) and the complexity of the expression (plots second row). Similar to \cite{lacava2021contemporarysymbolicregressionmethods}, the Wilcoxon signed-rank test was employed while using the Bonferroni correction for determining the threshold $\alpha$.  The color intensity corresponds to the degree of statistical significance, where grey indicates no significant difference between methods (i.e., method1 and method2 exhibit statistically equivalent performance). Darker shades represent stronger statistical significance.} 
  \label{plot5:significance}
\end{figure}

Regularization positively influences the complexity of the candidate solution by comparing the achieved scores between the different techniques, which gets even more significant when dimensional consistency dominates the search. Relating this to the results from Fig. \ref{plot1:r2-test}, we can emphasize that this loss in complexity negatively affects the predictive capabilities. Contrary to using the backpropagation to enforce homogeneity, it demonstrates lower complexities, here up to $45\%$ smaller than the standard method and up $22\%$ compared to the regularization. Moreover, the proposed technique exhibits better robustness in the presence of noise. The enhanced stability can be explained by the fact that noise primarily affects the data, while the underlying physical relationship of the feature, expressed through dimensional homogeneity, remains intact. \\\\
\noindent
\textbf{Statistical Significance}\\
The metrics compared in Fig. \ref{plot1:r2-test} and Fig. \ref{plot2:complexity} can be further evaluated with a focus on their statistical significance.  Fig. \ref{plot5:significance} illustrates heatmaps displaying the statistical comparisons across various noise levels ($\gamma \in \{0,0.01,0.1\}$), with the first row signifying the achieved $\mathbf{R}^2$-Score and the second row showing the results on the complexity of the approximated equations. Similar to \cite{lacava2021contemporarysymbolicregressionmethods}, we employed the Wilcoxon signed-rank test pairwise between the approaches (method1, method2), using the Bonferroni correction \cite{bonfer} to determine the threshold $\alpha$. The null hypothesis posits that method1 and method2 demonstrate equivalent performance, with rejection occurring when the p-value falls below $\alpha$. The values outlining the statistical significance for the $\mathbf{R}^2$-Score reflect the results presented in Fig. \ref{plot1:r2-test}. The application of varying penalty terms on the loss function for violating the homogeneity reveals marginal statistical significance $p<\alpha$ for $\gamma=0$ ($0.01d_\phi$ to $1d_\phi $ and $10_\phi$), with diminishing effects as target noise increases. In contrast, the proposed approach demonstrates strong significance to both regularization (${p<1.0e{-3}\cdot\alpha}$) and at least (${p<1.0e{-1}\cdot\alpha}$) to the GEP across the conducted noise levels. Examining the expression complexity at $\gamma=0$, regularization with the highest penalty term ($10d_\phi$) exhibits significantly different behavior from the standard method, while implementations with lower penalties ($1d_\phi$ and $0.1d_\phi$) yield results comparable to the standard approach. As target noise increases, the complexity diverges more. Notably, analogous to the $\mathbf{R}^2$-Score analysis, the semantic backpropagation method demonstrates statistically significant differences from the comparative method across all noise levels. 

\subsubsection{SRBench Comparision}
\begin{figure}[t]
\centering
  \hspace*{0.75cm}\includegraphics[scale=0.45]{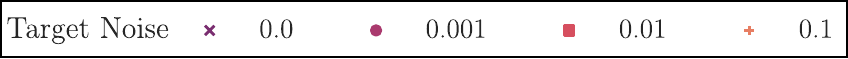}\\
  \vspace{0.2cm}
  \centering
  \includegraphics[scale=0.33]{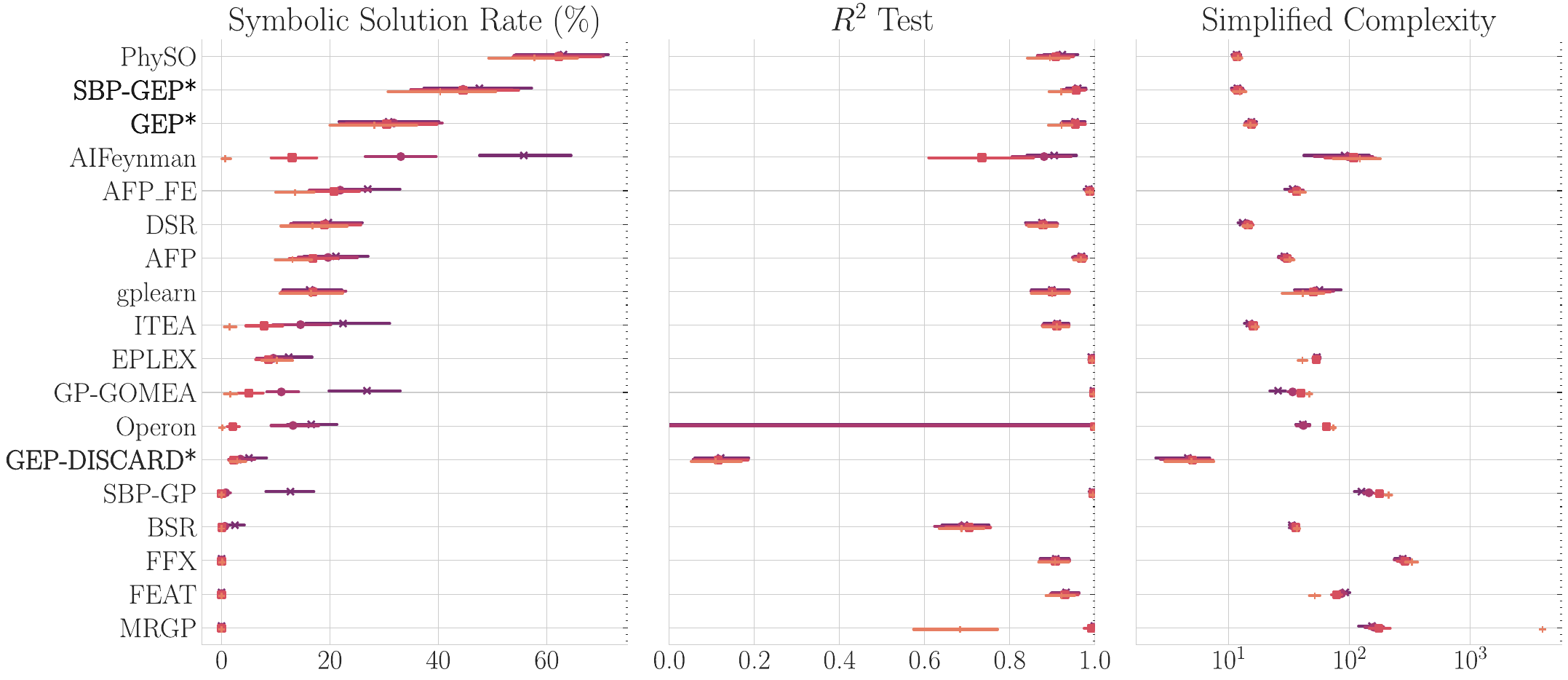}
  \caption{Pairplot comparing the performance of various methods listed in SRBench \cite{lacava2021contemporarysymbolicregressionmethods} and recently published \cite{Tenachi2023} utilizing the Feynman dataset with different levels of target noise ($\gamma \in \{0,0.001,0.01,0.1\}$). Metrics include symbolic solution rate, the $\textbf{R}^2$-test score, and the simplified complexity. SBP-GEP$^*$ denotes the GEP method utilizing semantic backpropagation to constrain solution with dimensional homogeneity, GEP$^*$ represents our basic implementation of Ferreira's \cite{ferreira2006gene} approach, and GEP-Discard$^*$ signifies a naive approach of discarding invalid candidates.}
  \label{plot3:metricSR}
\end{figure}

Fig. \ref{plot3:metricSR} summarizes the performance of the GEP and the proposed extension (outlined as `SBP-GEP') against the listed methods within the SRBench, and the results published by Tenachi et al. \cite{Tenachi2023}. Moreover, we included a method that utilizes the technique of discarding invalid candidate solutions (`GEP-Discard'). Here, the ranking (y-axis) is elucidated by the symbolic solution rate (Eqn. \ref{function:solacc}) in descending order. It further outlines the statistics according to the achieved $\mathbf{R}^2$ when employing the candidate solutions on the test set and the simplified complexity. Regarding the symbolic solution rate, it is noteworthy that the standard GEP implementation already outperforms the other approaches. The performance significantly deteriorates when employing a discarding approach (`GEP-Discard'), with a solution rate decreasing by 89\% and the $\mathbf{R}^2$ score dropping by 87\%. This degradation can be attributed to two main factors: (1) the absence of valid individuals in the initial population and (2) the limited number of valid candidate solutions during optimization. For case (1), the optimization stagnated with all individuals showing $\inf$ fitness, while case (2) suffered from a low diversity within the population that prevented a wide exploration.

When focusing specifically on the different evolutionary methods, this translates to a rate that, on average, surpasses the AFP\_FE by $45\%$, gplearn by $81\%$, the GP-GOMEA approximately more than two times, and the SBP-GP more than eight times. Akin to gplearn, the chosen base approach for our extension (GEP) demonstrates a low dispersion under inferring noise (from $\gamma=0$ to $\gamma=0.1$), where the solution rate faces a drop of $8.59\%$ (from $30.84$ - $28.19$), compared to AFP\_FE with $50\%$ or GP-GOMEA with $93\%$.

Using the information on dimensional homogeneity, we can further improve the average performance of the GEP by $33\%$, from $30.3$ to $44.27$. Referring to the other methods, only the PhySO approach, also incorporating the physical dimensions, demonstrates a higher solution rate (up to $38\%$). Apart from the fact that this approach shows a higher solution rate when measuring the $\textbf{R}^2$-test-score, the proposed method performs slightly better in this category. When focusing on the simplified complexity, the equations explored by the SBP-GEP range approximately at the same level. 

\section{Conclusion \& Future work}
This study proposed a methodology for incorporating constraints through semantic backpropagation into GEP. The example explored here is ensuring dimensional homogeneity. We adapted and included a vector representation for assigning properties to our input features, allowing the application of algebraic operations and introducing a distance metric. To determine the effectiveness of our method and the quality of the resultant solutions, a series of tests on ten different trials using the Feynman data set was conducted. We first compared the methodology against the most common strategy (regularization) and further against the methods listed in the SRBench.  

The results clearly outline the better performance compared to the regularization method, leading to a better score and more coherent candidate solutions regarding their functional complexity. Furthermore, results from the SRBench emphasize that the developed method and the implemented standard can outperform other existing approaches from the literature. 

Enhancing the GEP by introducing a dimension-aware semantic backpropagation mechanism enables the algorithm to systematically propagate information through the expression tree, facilitating targeted modifications that guide dimensional homogeneity. Therefore, this evolutionary method now generates approximations aligning with the physical dimensions, bridging the gap between data-driven exploration and physics-informed modeling.

Further work would be considered to improve the algorithm further. A critical aspect of the methodology is the existence of appropriate solutions stored in a library. Instead of random tree permutations, grammar-guided approaches \cite{crochepierre2022reinforcement} or language model-based approaches \cite{reissmann2024accelerating} could be applied. Furthermore, the operator can be extended and reformulated to satisfy further relations or symmetries of feature compositions. 
\backmatter
\bmhead{Acknowledgements:} This research was supported by The University of Melbourne’s Research Computing Services and the Petascale Campus Initiative. The work was supported by a Melbourne Research Scholarship provided by the University of Melbourne.

\section*{Declarations}

\subsection*{Conflicts of Interest}
The authors declare no conflicts of interest.

\subsection*{Code Availability}
Code and experimental setup available: \href{https://github.com/maxreiss123/GeneExpressionProgramming.jl}{GeneExpressionProgramming.jl}


\begin{appendices}
\label{section:appendix}
\newpage
\section{Entities of Gene Expression Programming}\label{secA}
The subsequent part introduces the key entities utilized within the GEP, which are motivated by leveraging advantages from the Genetic Algorithm (GA) \cite{holland-ga} and the Genetic Programming (GP) \cite{koza1992genetic}. While the first mentioned supports an efficient application of the genetic operators, the second one finds its strength in describing complex expression trees \cite{ferreira2006gene}. Therefore, this approach distinguishes between the genotype (genome) and the phenotype (expression tree).
\subsection{Genotype}
Akin to GA, each candidate solution within GEP is encoded as a linear string of fixed length but employing arbitrary symbols rather than a binary symbol set. This linear representation, the genome, can be further separated into genes, each consisting of a head and a tail section. Sequence \ref{first_karv} illustrates a typical gene structure:
\begin{align}
\label{first_karv}
{ \text{\tiny 0}\atop \sqrt()}\, {\text{\tiny 1} \atop *}\, {\text{\tiny 2} \atop +}\, {\text{\tiny 3}\atop a}\,{\text{\tiny 4} \atop b}\,{\text{\tiny 5} \atop -}\,{\text{\tiny 6} \atop c}\,{\text{\tiny 7} \atop d}\,{\text{\tiny 8} \atop a}\,{\text{\tiny 9} \atop b}\,{\text{\tiny 10} \atop c}\,{\text{\tiny 11} \atop d}\,{\text{\tiny 12} \atop d}
\end{align}
The gene structure comprises two distinct regions:
\begin{itemize}
\item Head [positions 0-5]: Contains both terminal and non-terminal symbols
\item Tail [positions 6-12]: Contains exclusively terminal symbols
\end{itemize}
To ensure a valid tree construction during decoding, the tail length must satisfy the relationship: $tail = head + 1$.
\subsection{Phenotype}
The decoding process constructs the expression tree by traversing the sequence \ref{first_karv} from left to right until all leaf nodes are populated with terminal symbols. Consequently, not all symbols in the gene sequence are necessarily incorporated into the final expression tree. For the example given \ref{sec_karv}, the decoding terminates at position seven, with the positions [0-7] forming the Karva expression. The unused symbols (Equation \ref{sec_karv} $I$) remain in the genome, serving as material in case genetic operators introduce additional non-terminal symbols into the head region [0-5].
\begin{align}
\label{sec_karv}
{ \text{\tiny 0}\atop \sqrt()}\, {\text{\tiny 1} \atop *}\, {\text{\tiny 2} \atop +}\, {\text{\tiny 3}\atop a}\,{\text{\tiny 4} \atop b}\,{\text{\tiny 5} \atop -}\,{\text{\tiny 6} \atop c}\,{\text{\tiny 7} \atop d}\,{\atop |}\underbrace{{\text{\tiny 8} \atop a}\,{\text{\tiny 9} \atop b}\,{\text{\tiny 10} \atop c}\,{\text{\tiny 11} \atop d}\,{\text{\tiny 12} \atop d}}_{(I)}
\end{align}
Based on the findings in \cite{Li2005PrefixGE}, the construction is performed in pre-order instead of the initially utilized level-order method, resulting in an expression tree illustrated in Fig. \ref{fig4:tree-appendix}:
 \begin{figure}[t]
  \centering
  \includegraphics[scale=.7]{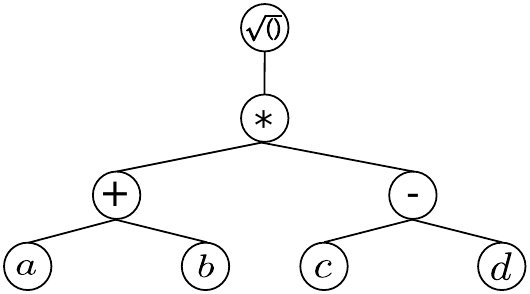}
  \caption{Expression tree built up in pre-order presents the phenotype for the sequence in Eqn. \ref{sec_karv}.
  }
  \label{fig4:tree-appendix}
\end{figure}

\section{Basics of Physical Dimension}\label{secC}
The following section presents a comprehensive overview of the International System of Units (SI).
\subsection{Base Units}
The system comprises seven units (or physical dimensions) constituting a foundation for physical measurement and calculations (Table \ref{tab:dimensions} \cite{Wiersma2021}):
\begin{table}[h]
\centering
\caption{Base Dimensions and Corresponding SI Units}
\begin{tabular}{|l|l|l|l|}
\hline
\textbf{Physical Quantity} & \textbf{Dimension} & \textbf{Unit Name} & \textbf{Symbol} \\
\hline
Mass & $[M]$ & kilogram & kg \\
Length & $[L]$ & meter & m \\
Time & $[T]$ & second & s \\
Temperature & $[\Theta]$ & kelvin & K \\
Amount of Substance & $[N]$ & mole & mol \\
Electric Current & $[I]$ & ampere & A \\
Luminous Intensity & $[J]$ & candela & cd \\
\hline
\end{tabular}
\label{tab:dimensions}
\end{table}
\\
Under applying certain mathematical functions, the rules as listed in Table \ref{tab:dim_operations} are utilized:
\begin{table}[h]
\centering
\caption{Rules for Dimensional Analysis Operations}
\label{tab:dim_operations}
\begin{tabular}{|l|l|l|}
\hline
\textbf{Operation} & \textbf{Rule} & \textbf{Example} \\
\hline
Addition/Subtraction & Same dimensions only & $[M] + [M] = [M]$ \\
Multiplication & Dimensions multiply & $[M] \cdot [LT^{-2}] = [MLT^{-2}]$ \\
Division & Dimensions divide & $[MLT^{-2}]/[L^2] = [ML^{-1}T^{-2}]$ \\
Powers & Exponents apply to all factors & $[ML]^2 = [M^2L^2]$ \\
Square Root & Halves all exponents & $\sqrt{[M^2L^2]} = [ML]$ \\
\hline
\end{tabular}
\end{table}

\subsection{Derived Units}
By employing these SI units and their corresponding mathematical operations, further units can be established \cite{taylor1995guide}. Table \ref{tab:derived_dimensions} highlights frequently encountered dimension compositions:

\begin{table}[htb]
\centering
\caption{Common Encountered Quantities: Units/ Dimensions}
\begin{tabular}{|l|l|l|l|l|}
\toprule
\textbf{Physical Quantity} & \textbf{Dimension} & \textbf{Unit Name} & \textbf{Symbol} & \textbf{SI Base Units} \\
\midrule
Force & $[ML T^{-2}]$ & newton & N & kg$\cdot$m$\cdot$s$^{-2}$ \\
Energy & $[M L^2 T^{-2}]$ & joule & J & kg$\cdot$m$^2\cdot$s$^{-2}$ \\
Power & $[M L^2 T^{-3}]$ & watt & W & kg$\cdot$m$^2\cdot$s$^{-3}$ \\
Pressure & $[M L^{-1} T^{-2}]$ & pascal & Pa & kg$\cdot$m$^{-1}\cdot$s$^{-2}$ \\
Velocity & $[L T^{-1}]$ & meter per second & m/s & m$\cdot$s$^{-1}$ \\
Acceleration & $[L T^{-2}]$ & meter per second squared & m/s$^2$ & m$\cdot$s$^{-2}$ \\
Density & $[M L^{-3}]$ & kilogram per cubic meter & kg/m$^3$ & kg$\cdot$m$^{-3}$ \\
Electric Charge & $[I T]$ & coulomb & C & A$\cdot$s \\
Electric Potential & $[M L^2 T^{-3} I^{-1}]$ & volt & V & kg$\cdot$m$^2\cdot$s$^{-3}\cdot$A$^{-1}$ \\
Resistance & $[M L^2 T^{-3} I^{-2}]$ & ohm & $\Omega$ & kg$\cdot$m$^2\cdot$s$^{-3}\cdot$A$^{-2}$ \\
Magnetic Flux & $[M L^2 T^{-2} I^{-1}]$ & weber & Wb & kg$\cdot$m$^2\cdot$s$^{-2}\cdot$A$^{-1}$ \\
\bottomrule
\end{tabular}
\label{tab:derived_dimensions}
\end{table}

\section{Semantic Backpropagation}\label{secD}
The subsequent contains supplementary material outlining technical details according to the proposed method.
\subsection{Backward Rules}
Akin to the training of neural networks, the backward pass is employed in many iterations to adjust the elements of the computational graph \cite{Wang2023}. We must define operators propagating information back to emulate this for the introduced method. Table \ref{tab:symsops} indicates that for operators like `$+$' or `$-$', this is straightforward because these do not transform intermediate states of $\phi$. In contrast, the multiplication and division operation has a significant impact. Therefore, their rules consist of defined algorithmic steps that consider various cases to determine which fraction of the target dimension should be further propagated to the corresponding underlying sub-tree. 
\begin{itemize}
\item Backward rule multiplication depicted in algorithm \ref{algo:backw-mul}:
\begin{algorithm}[H]
\caption{Perform Multiplication Rule}
\label{algo:backw-mul}
\begin{algorithmic}[1]
\Require $\hat{f}_i, \phi_c$ \Comment{$\phi_c$ target dimension} 
\Ensure $\phi_a, \phi_b$  \Comment{calculated residual target dimensions}
\State $\text{g}_i \Leftarrow \text{convert}(\hat{f}_i)$
\State $\phi_{tl} \Leftarrow \phi(get\_node\_left(\text{g}_i))$
\State $\phi_{tr} \Leftarrow 
\phi(get\_node\_right(\text{g}_i))$
\If{$d_\phi(\phi_{tl}, \phi_c) < \epsilon$} 
    \State $\phi_a \Leftarrow d_\phi(\phi_{tl}$
    \State $\phi_b \Leftarrow \phi_c - \phi_a$ \Comment{Employing the inverse to the forward operation}
\ElsIf{$d_\phi(\phi_{rl}, \phi_c) < \epsilon$}
    \State $\phi_b \Leftarrow d_\phi(\phi_{tl})$
    \State $\phi_a \Leftarrow \phi_c - \phi_b$ \Comment{Employing the inverse to the forward operation}
\Else \Comment{Split the dimension for further propagation}
    \State $\phi_a \Leftarrow \phi_c - \phi_c \div 2$
    \State $\phi_b \Leftarrow \phi_c - \phi_a$
\EndIf
\State \Return $\phi_a, \phi_b$
\end{algorithmic}
\end{algorithm}
\item Backward rule division outlined in algorithm  \ref{algo:backw-div}:
\begin{algorithm}[H]
\caption{Perform Division Rule}
\label{algo:backw-div}
\begin{algorithmic}[1]
\Require $\hat{f}_i, \phi_c$ \Comment{$\phi_c$ target dimension} 
\Ensure $\phi_a, \phi_b$  \Comment{calculated residual target dimensions}
\State $\text{g}_i \Leftarrow \text{convert}(\hat{f}_i)$
\State $\phi_{tl} \Leftarrow \phi(get\_node\_left(\text{g}_i))$
\State $\phi_{tr} \Leftarrow 
\phi(get\_node\_right(\text{g}_i))$
\If{$d_\phi(\phi_{tl}, \phi_c) < \epsilon$} 
    \State $\phi_a \Leftarrow d_\phi(\phi_{tl}$
    \State $\phi_b \Leftarrow -(\phi_c + \phi_a)$ \Comment{Employing the inverse to the forward operation}
\ElsIf{$d_\phi(\phi_{rl}, \phi_c) < \epsilon$}
    \State $\phi_b \Leftarrow d_\phi(\phi_{tl})$
    \State $\phi_a \Leftarrow \phi_c + \phi_b$ \Comment{Employing the inverse to the forward operation}
\Else \Comment{Split the dimension for further propagation}
    \State $\phi_a \Leftarrow \phi_c - \phi_c \div 2$
    \State $\phi_b \Leftarrow -(\phi_c + \phi_a)$
\EndIf
\State \Return $\phi_a, \phi_b$
\end{algorithmic}
\end{algorithm}
\end{itemize}

\subsection{Function Library}
Applying a search operator based on semantic backpropagation as proposed in \cite{Pawlak2015} requires a library to store the desired candidate solution. This library can be initialized in advance or filled up during runtime. Referring to a problem from the domain of symbolic regression, where we can describe such stored representations as trees, the question of the number of stored subsolutions, the corresponding memory overhead, the consistency, and the acceptable amount of time for fetching these intermediate states arises: 
\begin{figure}[t]
\centering
  \label{fig:timelib}
  \hspace*{0.75cm}\includegraphics[scale=0.45]{images/pairplot_symbolic_legend.pdf}\\
  \vspace{0.2cm}
  \centering
  \includegraphics[scale=0.33]{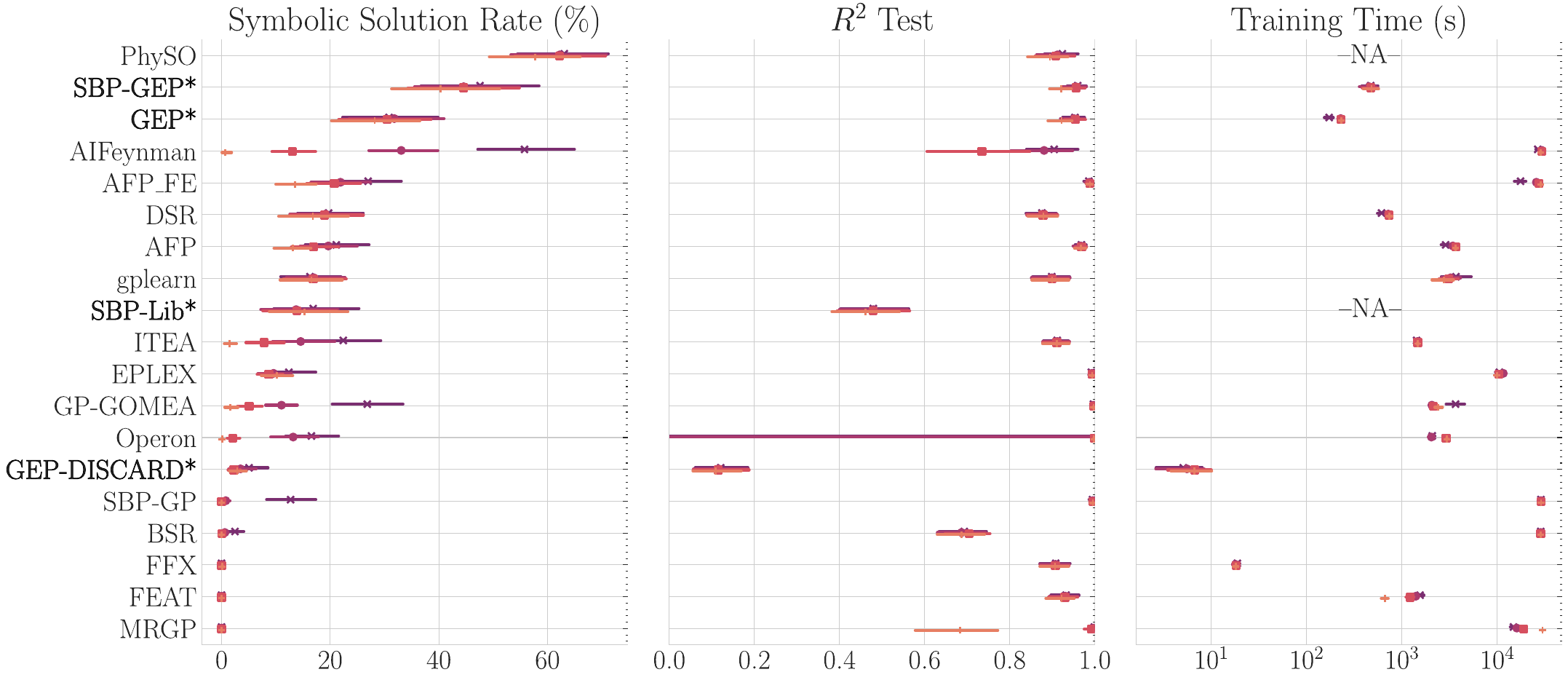}
  \caption{Pairplot comparing the performance of various methods listed in SRBench \cite{lacava2021contemporarysymbolicregressionmethods} and recently published \cite{Tenachi2023} utilizing the Feynman dataset on different levels of target noise ($\gamma \in \{0,0.001,0.01,0.1\}$). Metrics include the symbolic solution rate, the $\textbf{R}^2$-test-score, and the execution time. The latter metric was multiplied by the threads employed. SBP-GEP$^*$ denotes the GEP method utilizing semantic backpropagation to ensure dimensional homogeneity, SBP\_Lib$^*$ the dimension-based library, while GEP$^*$ represents our basic implementation of Ferreira's \cite{ferreira2006gene} approach.}
  \label{plot:metricSR-appendix}
\end{figure}
\begin{enumerate}
\item \textbf{Memory overhead}: The overhead of the memory once depends on the total number of trees and the internal tree representation. In our GEP implementation, we have chosen a more data-oriented design. Instead of focusing on code readability or object abstraction, the gene information is stored as an integer representation in an array (max alignment eight bytes per entry), whereby the size is restricted by the head length ($hl$) (hyperparameter \ref{tab:high}). Constructing the library is an iterative process, inserting new non-terminals and complementing them with terminals or the sub-trees created. Leveraging this process, an upper limit can be estimated by conducting the theoretical analysis of search space scaling outlined by Jiang et al. \cite{10.1007/978-3-031-43421-1_11}. Given the terminals as $m$ and the non-terminals as $o$, the size of the search space with $l$ internal nodes scales with:
\begin{align}
    \mathcal{O}\left((4(m+1)o)^{\frac{l-1}{2}} \right) \,.
\end{align}
Based on that, the upper limit can be estimated:
\begin{align}
    |Lib(o,m,hl)| = \sum_i^{hl} \left( \min(10^5, (4(m+1)o)^{\frac{hl-1}{2}})\right) 
\end{align}
where $10^5$ limits the amount of newly generated candidates according to practical reasons. Assuming $m=7$ as the average inputs, $hl=8$ and $o=9$ lead to approximately $3.9*10^5$. Moreover, the array of length $hl$ and the mentioned alignment (8 Byte per position using Int64) yields a library of around 25Mb memory. Which can be further optimized by considering UInt32 or UInt8. \\
\item \textbf{Overhead Computational Time}:
Applying a correction mechanism to a large proportion of the population raises further questions about the effect of the runtime, including the time for the library creation and execution. The overhead was quantified by monitoring the execution throughout the experimental process. Fig. (\ref{fig:timelib}) (right-panel) provides an indication of the required time, including the different noise levels multiplied by the threads utilized. The results demonstrate while the proposed method offers superior performance compared to the standard approach, it extends the execution by approximately two to 2.5 times. 
\\
\item \textbf{Consistency}: Library construction adheres to predefined rules to avoid homogeneity violations (ref. Eqn. \ref{function:dimension}), preventing arithmetic issues or the occurrence of redundant nested functions like ($\exp(\log(\cos(...)))$). Since these entries are utilized to accomplish a beneficial exploration, we expect them to be appropriate for the given task. To test this hypothesis, we conduct the same test setup as in \ref{results_and_dists} and compare the $\mathbf{R}^2$ (Eqn. \ref{eqn:r2-score}) as well the solution rate (Eqn. \ref{function:solacc}). In detail, we employ all entries that match the target physical dimension in combination with coefficient optimization. The result presented in Fig. \ref{plot:metricSR-appendix} reveals that the entries outer show similar performance to the gplearn. Relating this to the amount of complexity that can show up and the number of equations within the `easy' set implies that the content can at least discover approximately $50\%$ from the mentioned subset. On the contrary, the $\mathbf{R}^2$ is dominated by the low accuracy of the candidate solutions for `medium' and `hard' cases.

\end{enumerate}
\end{appendices}
\newpage
\bibliography{bib/main}
\end{document}